\documentclass[11pt]{article}
\usepackage{amsfonts}
\usepackage{mathrsfs}
\usepackage[dvips]{graphics}
\usepackage[dvips]{color}
\usepackage{amsthm}
\usepackage[]{amsmath}
\usepackage{CJK}
\usepackage{indentfirst}
\newtheorem{proposition}{Proposition}[section]

\newtheorem{theorem}[proposition]{Theorem}

\newtheorem{property}[proposition]{Property}
\newtheorem{remark}[proposition]{Remark}

\begin{document}
\begin{CJK*}{GBK}{song}
\CJKindent

\centerline{\textbf{\LARGE{The Numbers of Distinct Squares}}}

\vspace{0.2cm}

\centerline{\textbf{\LARGE{and Cubes in the Tribonacci Sequence}}}

\vspace{0.2cm}

\centerline{Huang Yuke\footnote[1]{School of Mathematics and Systems Science, Beihang University (BUAA), Beijing, 100191, P. R. China. E-mail address: huangyuke07@tsinghua.org.cn,~hyg03ster@163.com (Corresponding author).}
~~Wen Zhiying\footnote[2]{Department of Mathematical Sciences, Tsinghua University, Beijing, 100084, P. R. China. E-mail address: wenzy@tsinghua.edu.cn.}}

\vspace{1cm}

\noindent\textbf{Abstract:} The Tribonacci sequence $\mathbb{T}$ is the fixed point of the substitution $\sigma(a,b,c)=(ab,ac,a)$.
In this note, we give the explicit expressions of the numbers of distinct squares and cubes in $\mathbb{T}[1,n]$ (the prefix of $\mathbb{T}$ of length $n$).

\vspace{0.2cm}


\vspace{0.2cm}

\noindent\textbf{2010 MR Subject Classification:} 11B85; 68Q45

\section{Introduction}

Let $\mathcal{A}=\{a,b,c\}$ be a three-letter alphabet.
The concatenation of factors $\nu$ and $\omega$ denoted by $\nu\ast\omega$ or $\nu\omega$.
The Tribonacci sequence $\mathbb{T}$ is the fixed point beginning with $a$ of the substitution $\sigma$ defined by $\sigma(a)=ab$, $\sigma(b)=ac$, $\sigma(c)=a$.
The Tribonacci sequence, which is a natural generalization of the Fibonacci sequence, has
been studied extensively by many authors, see \cite{DM2008,MS2014,RT2005,RSZ2010,TW2007}.
%AR1991,CHM2001,

Let $\omega$ be a factor of $\mathbb{T}$, denoted by $\omega\prec\mathbb{T}$.
Let $\omega_p$ be the $p$-th occurrence of $\omega$.
If the factor $\omega$ and integer $p$ such that $\omega_p\omega_{p+1}$ (resp. $\omega_p\omega_{p+1}\omega_{p+2}$) is the factor of $\mathbb{T}$, we call it a square (resp. cube) of $\mathbb{T}$.
As we know, $\mathbb{T}$ contains no fourth powers. The properties of squares and cubes are objects of a great interest in many aspects of mathematics and computer science etc.

We denote by $|\omega|$ the length of $\omega$.
Let $\tau=x_1\cdots x_n$ be a finite word (or $\tau=x_1x_2\cdots$ be a sequence).
For any $i\leq j\leq n$, define $\tau[i,j]:=x_ix_{i+1}\cdots x_{j-1}x_j$.
By convention, denote $\tau[i]:=\tau[i,i]=x_i$, $\tau[i,i-1]:=\varepsilon$(empty word).
Denote $T_m=\sigma^m(a)$ for $m\geq0$, $T_{-2}=\varepsilon$, $T_{-1}=c$,
then $T_0=a$, $T_1=ab$ and $T_{m}=T_{m-1}T_{m-2}T_{m-3}$ for $m\geq2$.
Denote $t_m=|T_m|$ for $m\geq-2$, called the $m$-th Tribonacci number.
Then $t_{-2}=0$, $t_{-1}=t_{0}=1$, $t_m=t_{m-1}+t_{m-2}+t_{m-3}$ for $m\geq1$.
Denote by $\delta_m$ the last letter of $T_m$, then $\delta_m=a$ (resp. $b$, $c$) for $m\equiv0$, (resp. $1$, $2$) mod 3, $m\geq-1$.

In 2006, A.Glen gave the number of distinct squares in $T_m$ in her PhD thesis, see Theorem 6.30 in \cite{G2006}.
She also gave some properties about cubes.
In 2014, H.Mousavi and J.Shallit\cite{MS2014} gave some properties of the Tribonacci word, such as the lengths of squares and cubes.
All of these results above only consider the squares or cubes in the prefixes of some special lengths: the Tribonacci numbers.
The main aim of this article is to give the explicit expressions of the numbers of distinct squares and cubes in $\mathbb{T}[1,n]$.
This topic has been studied not only for the Tribonacci sequence.
We gave the expressions of the numbers of distinct squares and cubes in each prefix of the Fibonacci sequence, see \cite{HW2016-2}.

The main tools of the paper are ``kernel word'' and the unique decomposition
of each factor with respect to its kernel, which introduced and studied in Huang and Wen\cite{HW2015-2}.
In this paper, Section 2 presents some basic notations and known results.
Section 3 proves some basic properties of squares.
We determine the numbers of distinct squares and cubes in $\mathbb{T}[1,n]$ in
Section 4 and 5 respectively.
In Section 6, we give some open problems.

\section{Preliminaries}

We define the kernel numbers that $k_{0}=0$, $k_{1}=k_{2}=1$, $k_m=k_{m-1}+k_{m-2}+k_{m-3}-1$ for $m\geq3$.
The kernel word with order $m$ is defined as
$K_1=a$, $K_2=b$, $K_3=c$, $K_m=\delta_{m-1}T_{m-3}[1,k_m-1]\text{ for }m\geq4.$
By Proposition 2.7 in \cite{HW2015-2}, all kernel words are palindromes.
Let $Ker(\omega)$ be the maximal kernel word occurring in factor $\omega$, then by Theorem 4.3 in \cite{HW2015-2}, $Ker(\omega)$ occurs in $\omega$ only once.

\begin{property}[Theorem 4.11 in \cite{HW2015-2}]\label{wp}~
$\forall\omega\in\mathbb{T}$, $p\geq1$, $Ker(\omega_p)=Ker(\omega)_p$.
\end{property}

This means, let $Ker(\omega)=K_m$, then the maximal kernel word occurring in $\omega_p$ is just $K_{m,p}$. For instance, $Ker(aba)=b$, $(aba)_3=\mathbb{T}[8,10]$, $(b)_3=\mathbb{T}[9]$, so $Ker((aba)_3)=(b)_3$, $(aba)_3=a(b)_3a$.

The next three properties can be proved easily by induction.

\begin{property}[]\label{bp1}~For $m\geq3$,
(1) $k_m=k_{m-3}+t_{m-4}=\frac{t_{m-3}+t_{m-5}+1}{2}$;
%=k_{m-1}+t_{m-5}

(2) $K_m=\delta_{m-1}T_{m-4}K_{m-3}[2,k_{m-3}]=\delta_{m-1}T_{m-4}T_{m-5}[1,k_{m-3}-1]$.
\end{property}

\begin{property}[]\label{bp3}~(1) $\sum_{i=0}^{m}t_i=\frac{t_m+t_{m+2}-3}{2}$
for $m\geq0$, (see Lemma 6.7 in \cite{G2006});

(2) $\sum_{i=1}^{m}k_i=\frac{k_m+k_{m+2}+m-1}{2}=\frac{t_{m-2}+t_{m-3}+m}{2}$ for $m\geq1$.
\end{property}

\begin{property}[]\label{bp2}~For $m\geq0$,
(1) $T_{m}T_{m+1}[1,k_{m+4}-2]=T_{m+2}[1,k_{m+5}-2]$;

(2) $T_{m+3}[1,k_{m+6}-2]=T_{m+1}T_{m}T_{m+1}[1,k_{m+4}-2]$.
\end{property}

\section{Basic properties of squares}

We denote the gap between $\omega_p$ and $\omega_{p+1}$ by $G_p(\omega)$.
The sequence $\{G_p(\omega)\}_{p\ge 1}$ is called the gap sequence of factor $\omega$.
By Lemma 4.7, Definition 4.12 and Corollary 4.13 in \cite{HW2015-2}, any factor $\omega$ with kernel $K_m$ can be expressed uniquely as
$$\omega=T_{m-1}[i,t_{m-1}-1]K_mT_{m}[k_m,k_m+j-1],$$
where $1\leq i\leq t_{m-1}$ and $0\leq j\leq t_{m-1}-1$.
By Theorem 3.3, Corollary 3.12 and Proposition 6.7(1) in \cite{HW2015-2}, $\omega_p\omega_{p+1}\prec\mathbb{T}$ has three cases.

\vspace{0.2cm}

\textbf{Case 1.} $G_p(K_m)=G_1(K_m)$ and $|G_1(K_m)|=t_{m}-k_{m}$.

Since $|G_p(K_m)|=|T_{m-1}[i,t_{m-1}-1]|+|T_{m}[k_m,k_m+j-1]|=t_{m-1}-i+j$, $j=t_{m-2}+t_{m-3}-k_{m}+i$. So $0\leq j\leq t_{m-1}-1$ gives a range of $i$. Comparing this range with $1\leq i\leq t_{m-1}$, we have $1\leq i\leq k_{m+1}-1$ and $m\geq3$.
Furthermore,
\begin{equation*}
\begin{split}
\omega=&T_{m-1}[i,t_{m-1}]T_{m}[1,t_{m-2}+t_{m-3}+i-1]\\
=&T_{m-1}[i,t_{m-1}]T_{m-2}T_{m-3}T_{m-2}[1,i-1];\\
\omega\omega
=&T_{m}[i,t_{m}-1]\underline{\delta_{m}T_{m-1}[1,k_{m+1}-1]}T_{m+1}[k_{m+1},t_{m}+i-1].
\end{split}
\end{equation*}
Thus $K_{m+1}=\delta_{m}T_{m-1}[1,k_{m+1}-1]\prec\omega\omega$.
Similarly, $K_{m+2},K_{m+3},K_{m+4}\!\not\prec\omega\omega$ and $|\omega\omega|<K_{m+5}$, so $K_{m+1}$ is the largest kernel word in $\omega\omega$, i.e. $Ker(\omega\omega)=K_{m+1}$ for $m\geq3$.
Moreover, since $|\omega|=|G_p(K_m)|+k_m$, $|\omega|=t_{m}$.

\vspace{0.2cm}

By analogous arguments, we have

\textbf{Case 2.} $G_p(K_m)=G_2(K_m)$ and $|G_2(K_m)|=t_{m-2}+t_{m-1}-k_{m}$.
$$\omega\omega
=T_{m}[i,t_{m-1}+t_{m-2}-1]\underline{K_{m+2}}T_{m+1}[k_{m+2},t_{m-1}+t_{m-2}+i-1],$$
where $1\leq i\leq k_{m+2}-1$, $m\geq2$, $Ker(\omega\omega)=K_{m+2}$ and $|\omega|=t_{m-2}+t_{m-1}$.

\textbf{Case 3.} $G_p(K_m)=G_4(K_m)$ and $|G_4(K_m)|=t_{m-1}-k_{m}$.
$$\omega\omega
=T_{m-1}[i,t_{m-1}-1]\underline{K_{m+3}}T_{m+1}[k_{m+3},t_{m-1}+i-1],$$
where $k_{m}\leq i\leq t_{m-1}$, $m\geq1$, $Ker(\omega\omega)=K_{m+3}$ and $|\omega|=t_{m-1}$.

\begin{remark}
By the three cases of squares, we have: (1) all squares in $\mathbb{T}$ are of length $2t_m$ or $2t_m+2t_{m-1}$ for some $m\geq0$; (2) for all $m\geq0$, there exists a square of length $2t_m$ and $2t_m+2t_{m-1}$ in $\mathbb{T}$. These are known results of H.Mousavi and J.Shallit, see Theorem 5 in \cite{MS2014}.
\end{remark}

Denote $P(\omega,1)$ (resp. $L(\omega,1)$) the position of the last (resp. first) letter of $\omega_1$.
By Theorem 6.1, Remark 6.2 in \cite{HW2015-2} and $P(\omega,1)=L(\omega,1)+|\omega|-1$, we have

\begin{property}\label{P}
$P(K_m,1)=t_{m-1}+k_m-1=k_{m+3}-1$ for $m\geq1$.
\end{property}

We define three sets for $m\geq4$,
\begin{equation*}
\begin{cases}
\langle1,K_m\rangle:=
\{P(\omega\omega,1):Ker(\omega\omega)=K_m,|\omega|=t_{m-1},\omega\omega\prec\mathbb{T}\}\\
\langle2,K_m\rangle:=
\{P(\omega\omega,1):Ker(\omega\omega)=K_m,|\omega|=t_{m-4}+t_{m-3},\omega\omega\prec\mathbb{T}\}\\
\langle3,K_m\rangle:=
\{P(\omega\omega,1):Ker(\omega\omega)=K_m,|\omega|=t_{m-4},\omega\omega\prec\mathbb{T}\}
\end{cases}
\end{equation*}
Obviously these sets correspond the positions $P(\omega\omega,1)$ for the three cases of squares respectively.
By Property \ref{P},
$\langle1,K_m\rangle$ is equal to
\begin{equation*}
\begin{split}
&\{P(\omega,1):
\omega=T_{m-1}[i,t_{m-1}-1]K_{m}T_{m}[k_{m},t_{m-1}+i-1],1\leq i\leq k_{m}-1\}\\
=&\{P(K_m,1)+t_{m-1}-k_{m}+i,1\leq i\leq k_{m}-1\}
=\{2t_{m-1},\cdots,k_{m+4}-2\}.\\
\end{split}
\end{equation*}
Moreover
$\sharp\langle 1,K_m\rangle=\sharp\{1\leq i\leq k_{m}-1\}=k_{m}-1$. Similarly

\begin{property}[]\label{K1} For $m\geq4$,
\begin{equation*}
\begin{cases}
\langle1,K_m\rangle=\{2t_{m-1},\cdots,k_{m+4}-2\};\\
\langle2,K_m\rangle=\{2t_{m-1}-t_{m-2},\cdots,t_{m-1}+k_{m+2}-2\};\\
\langle3,K_m\rangle=\{k_{m+3}-1,\cdots,t_{m-1}+2t_{m-4}-1\}.
\end{cases}
\end{equation*}
Moreover $\sharp\langle 1,K_m\rangle=\sharp\langle 2,K_m\rangle=k_{m}-1$,
$\sharp\langle 3,K_m\rangle=t_{m-4}-k_{m-3}+1$.
\end{property}

For $m\geq3$, denote
$$\begin{array}{c}
\Delta_m:=\sum\limits_{i=4}^m\sharp\langle 1,K_i\rangle=\sum\limits_{i=4}^m\sharp\langle 2,K_i\rangle,
~\Theta_m:=\sum\limits_{i=4}^m\sharp\langle 3,K_i\rangle.
\end{array}$$
Obviously $\Delta_3=\Theta_3=0$. For $m\geq4$,
by $\sum_{i=0}^{m}t_i=\frac{t_m+t_{m+2}-3}{2}$ and
$\sum_{i=1}^{m}k_i=\frac{t_{m-2}+t_{m-3}+m}{2}$, we have
$\Delta_{m}=\frac{t_{m-2}+t_{m-3}-m}{2}$ and $\Theta_m=\frac{t_{m-2}-t_{m-3}+2t_{m-4}+m-6}{2}$.

\section{The number of distinct squares in $\mathbb{T}[1,n]$}

By Property \ref{K1}, $\langle i,K_m\rangle$ are pairwise disjoint, and each set contains some consecutive integers. Therefore we get a chain
$$\langle3,K_4\rangle,\langle2,K_4\rangle,\langle1,K_4\rangle,\langle3,K_5\rangle,\cdots,
\langle3,K_m\rangle,\langle2,K_m\rangle,\langle1,K_m\rangle,\cdots$$

Denote $a(n):=\sharp\{\omega:\omega\omega\prec\mathbb{T}[1,n],
\omega\omega\not\!\prec\mathbb{T}[1,n-1]\}$, then $a(n)=1$ iff $n\in\cup_{m\geq4}(\langle3,K_m\rangle\cup\langle2,K_m\rangle\cup\langle1,K_m\rangle)$. The ``$\cup$'' in this note always means disjoint union.
Since $\max\langle1,K_m\rangle+1=\min\langle3,K_{m+1}\rangle$, sets $\langle1,K_m\rangle$ and $\langle3,K_{m+1}\rangle$ are consecutive.
We have
$\langle1,K_m\rangle\cup\langle3,K_{m+1}\rangle=\{2t_{m-1},\cdots,t_{m}+2t_{m-3}-1\}$.

\begin{property}[]\label{a}~For $n<14$, $a(n)=1$ iff $n\in\{8,10\}$;
for $n\geq14$, let $m$ such that $2t_{m-1}\leq n<2t_{m}$, then $m\geq4$ and $a(n)=1$ iff
$$n\in\{2t_{m-1},\cdots,t_{m}+2t_{m-3}-1\}
\cup\{2t_{m}-t_{m-1},\cdots,t_{m}+k_{m+3}-2\}.$$
\end{property}

Denote $A(n):=\sharp\{\omega:\omega\omega\prec\mathbb{T}[1,n]\}$, which is the number of distinct squares in $\mathbb{T}[1,n]$. Then $A(n)=\sum_{i=1}^n a(i)$.
For $n\geq14$, find the $m$ such that $2t_{m-1}\leq n<2t_{m}$, then $m\geq4$.
Denote
\begin{equation*}
\begin{cases}
\alpha_m:=\min\langle1,K_m\rangle=2t_{m-1},~\beta_m:=\max\langle3,K_{m+1}\rangle=t_{m}+2t_{m-3}-1,\\
\gamma_m:=\min\langle2,K_{m+1}\rangle=2t_{m}-t_{m-1},\\
\theta_m:=\max\langle2,K_{m+1}\rangle=t_{m}+k_{m+3}-2=\frac{3t_{m}+t_{m-2}-3}{2}.
\end{cases}
\end{equation*}

By Property \ref{a} and the definition of $\Delta_m$, $\Theta_m$, for $m\geq4$
\begin{equation*}
\begin{cases}
A(\alpha_m)=\Delta_{m-1}+\Delta_{m}+\Theta_{m}+1
=\frac{2t_{m-2}+t_{m-3}+3t_{m-4}-m-3}{2},\\
A(\beta_m)=\Delta_{m}+\Delta_{m}+\Theta_{m+1}
=\frac{t_{m-1}+t_{m-2}+4t_{m-3}-m-5}{2},\\
A(\gamma_m)=\Delta_{m}+\Delta_{m}+\Theta_{m+1}+1
=A(\beta_m)+1,\\
A(\theta_m)=\Delta_{m}+\Delta_{m+1}+\Theta_{m+1}
=A(\alpha_{m+1})-1.
\end{cases}
\end{equation*}

Obviously, when $\alpha_m\leq n<\beta_m$, $A(n)=A(\alpha_m)+n-\alpha_m$;
when $\beta_m\leq n<\gamma_m$, $A(n)=A(\beta_m)$;
when $\gamma_m\leq n<\theta_m$, $A(n)=A(\gamma_m)+n-\gamma_m$;
when $\theta_m\leq n<\alpha_{m+1}$, $A(n)=A(\theta_m)$. So

\begin{theorem}[The numbers of distinct squares, $A(n)$]\label{A}\
$A(n)=0$ for $n\leq7$; $A(n)=1$ for $n=8,9$; $A(n)=2$ for $10\leq n\leq13$.
For $n\geq14$, let $m$ such that $\alpha_m\leq n<\alpha_{m+1}$, then $m\geq4$,
\begin{equation*}
A(n)=
\begin{cases}
n-\frac{1}{2}(t_{m}+t_{m-3}+m+3),&\alpha_m\leq n<\beta_{m};\\
\frac{1}{2}(t_{m-1}+t_{m-2}+4t_{m-3}-m-5),&\beta_{m}\leq n<\gamma_m;\\
n-\frac{1}{2}(t_{m-1}+3t_{m-2}+m+3),&\gamma_m\leq n<\theta_m;\\
\frac{1}{2}(2t_{m-1}+t_{m-2}+3t_{m-3}-m-6),&\theta_m\leq n<\alpha_{m+1}.
\end{cases}
\end{equation*}
\end{theorem}

\noindent\textbf{Example.} Count the numbers of distinct squares in $\mathbb{T}[1,355]$,
i.e. $A(355)$.

Since $n=355\geq14$ and $2t_{8}=298\leq n=355<2t_{9}=548$, $m=9$.
Moreover $n<\beta_9=t_{9}+2t_{6}-1=361$, so $C(355)=355-\frac{1}{2}(t_{9}+t_{6}+12)=190$.

\vspace{0.2cm}

For $m\geq3$, since $\theta_{m-1}=\frac{3t_{m-1}+t_{m-3}-1}{2}\leq t_m<\alpha_{m}$, by Theorem \ref{A}, we have $A(t_m)=A(\theta_{m-1})$ for $m\geq5$.
It is easy to check the expression holds also for $m=3,4$. Thus for $m=0,1,2$, $A(t_m)=0$, and

\begin{theorem}[]\ For $m\geq3$,
$A(t_m)=\frac{1}{2}(2t_{m-2}+t_{m-3}+3t_{m-4}-m-5)$.
\end{theorem}

\begin{remark} A.Glen show the number of distinct squares in $T_m$, i.e. $A(t_m)$ in Theorem 6.30 in \cite{G2006}, that
$A(t_m)=\sum_{i=0}^{m-2}(d_i+1)+d_{m-4}+d_{m-5}+1$ for $m\geq3$,
where $d_{-2}=d_{-1}=-1$, $d_0=0$ and $d_m=\frac{t_{m+1}+t_{m-1}-3}{2}$ for $m\geq1$.
By the expression of $\sum_{i=0}^{m}t_i$, we know the two expressions of $A(t_m)$ are same.
\end{remark}

\section{The number of distinct cubes in $\mathbb{T}[1,n]$}

Let $\omega$ be a factor with kernel $K_m$, by analogous arguments as Section 3 and 4,
we have: (1) all cubes in $\mathbb{T}$ are of length $3t_m$ for some $m\geq3$; (2) for all $m\geq3$, there exists a cube of length $3t_m$, which is Theorem 7 in \cite{MS2014}.
Moreover,

\begin{theorem}[The numbers of distinct cubes, $B(n)$]\label{B}\
$B(n)=0$ for $n\leq57$.
For $n\geq58$, let $m$ such that $t_{m-1}+2t_{m-4}\leq n<t_{m}+2t_{m-3}$, then $m\geq7$,
\begin{equation*}
B(n)=
\begin{cases}
n-\frac{1}{2}(4t_{m-1}-t_{m-2}-3t_{m-3}+m-6),&n\leq\frac{3t_{m-1}-t_{m-3}-3}{2};\\
\frac{1}{2}(t_{m-5}+t_{m-6}-m+3),&otherwise.
\end{cases}
\end{equation*}
\end{theorem}

\noindent\textbf{Example.} Count the numbers of distinct cubes in $\mathbb{T}[1,365]$,
i.e. $B(365)$.

Since $n=365\geq58$ and
$t_{9}+2t_{6}=362\leq n=365<t_{10}+2t_{7}=666$, $m=10$.
Moreover $n<t_{9}+k_{11}-2=369$,
$B(365)=365-\frac{4t_{9}-t_{8}-3t_{7}+4}{2}=11$.

\vspace{0.2cm}

For $m\geq7$, since $t_m>\frac{3t_{m-1}-t_{m-3}-3}{2}$, by Theorem \ref{B}, we have

\begin{theorem}[]\ For $m\leq6$, $B(t_m)=0$, for $m\geq7$,
$B(t_m)=\frac{t_{m-5}+t_{m-6}-m+3}{2}$.
\end{theorem}

\section{Open Problems}

In 2014, H.Mousavi and J.Shallit\cite{MS2014} gave explicit expressions about the numbers of repeated squares and cubes in the Tribonacci word $T_m$, which they proved by mechanical way.
In \cite{HW2016-2}, we give fast algorithms for counting the numbers of repeated squares and cubes in each prefix of the Fibonacci sequence.
But we have not yet succeeded in giving fast algorithms for counting the numbers of repeated
squares and cubes in $\mathbb{T}[1,n]$ for all $n$ by our method.

\end{CJK*}
\end{document}